\documentclass{amsart}

\usepackage{amssymb,latexsym,amsmath,amsfonts,hhline,comment,mathrsfs,stmaryrd}
\usepackage{pdfsync}
\usepackage{amsthm}
\usepackage{graphicx}
\usepackage{wrapfig}
\usepackage{caption}
\usepackage{subcaption}
\usepackage{indentfirst}
\usepackage{enumerate}
\input xy
\xyoption{all}

\newcommand{\cal}{\mathcal}

\usepackage{ifpdf}
\ifpdf \usepackage[colorlinks=true, citecolor=blue, linkcolor=blue, urlcolor=blue]{hyperref} \fi

\newtheorem{formula}{}[section]
\newtheorem{definition}[formula]{Definition}
\newtheorem{corollary}[formula]{Corollary}
\newtheorem{remark}[formula]{Remark}
\newtheorem{lemma}[formula]{Lemma}
\newtheorem{theorem}[formula]{Theorem}

\def\thrm{\begin{theorem}}
\def\thrml#1{\begin{theorem}\label{#1}}
\def\ethrm{\end{theorem}}
\def\rmrk{\begin{remark}}
\def\rmrkl#1{\begin{remark}\label{#1}}
\def\ermrk{\end{remark}}
\def\dfntn{\begin{definition}}
\def\dfntnl#1{\begin{definition}\label{#1}}
\def\edfntn{\end{definition}}
\def\nmrt{\begin{enumerate}}
\def\enmrt{\end{enumerate}}
\def\tm#1{\item[{\rm (#1)}]}
\def\qtnl#1{\begin{equation}\label{#1}}
\def\eqtn{\end{equation}}
\def\lmm{\begin{lemma}}
\def\lmml#1{\begin{lemma}\label{#1}}
\def\elmm{\end{lemma}}
\def\crllr{\begin{corollary}}
\def\crllrl#1{\begin{corollary}\label{#1}}
\def\ecrllr{\end{corollary}}
\def\css{\begin{cases}}
\def\ecss{\end{cases}}

\def\proof{\noindent{\bf Proof}.\ }

\def\cP{{\cal P}}
\def\cQ{{\cal Q}}
\def\cR{{\cal R}}

\def\cX{{\cal X}}
\def\cY{{\cal Y}}

\DeclareMathOperator{\dom}{Dom}

\DeclareMathOperator{\id}{id}

\DeclareMathOperator{\pr}{pr}

\def\eprf{\hfill$\square$}

\def\bone{{\bf 1}}
\def\grp#1{\langle {#1}\rangle}

\def\qaq{\quad\text{and}\quad}

\def\ov{\overline}

\def\VRT#1{*=<5mm>[o][F-]{#1}}
\def\VRTB#1{*=<9mm>[o][F-]{#1}}

\def\grphp#1{$\xymatrix@R=10pt@C=10pt@M=0pt@L=2pt{#1}$}

\begin{document}

\title{Tensor products of coherent configurations}
\author{Gang Chen}, 
\address{School of Mathematics and Statistics, Central China Normal University, Wuhan 430079, China}
\email{chengangmath@mail.ccnu.edu.cn}
\thanks{The first author is supported by NSFC grant No. 11971189}
\author{Ilia Ponomarenko}
\address{Steklov Institute of Mathematics at St. Petersburg, Russia;\newline Sobolev Institute of Mathematics, Novosibirsk, Russia; and\newline School of Mathematics and Statistics of Central China Normal University, Wuhan, China} 
\email{inp@pdmi.ras.ru}
\thanks{}
\date{}

\begin{abstract}
A Cartesian decomposition of a coherent configuration~$\cX$ is defined as a special set of its parabolics that form a Cartesian decomposition of the underlying set. It turns out that every tensor decomposition of~$\cX$ comes from a certain Cartesian decomposition. It is proved that if the coherent configuration $\cX$ is thick, then there is a unique  maximal  Cartesian decomposition of~$\cX$, i.e., there is exactly one internal tensor decomposition of~$\cX$ into indecomposable components.  In particular, this implies an analog of the Krull--Schmidt theorem for the thick coherent configurations. A polynomial-time algorithm for finding the maximal  Cartesian decomposition of a thick coherent configuration is constructed.
\end{abstract}

\maketitle

\section{Introduction}

A natural problem arising in the study of direct products in various categories is to find an analog of the Krull--Schmidt theorem which asserts  that a group subject to certain conditions can be essentially uniquely written as a finite direct product of indecomposable groups. Despite a general solution proposed in~\cite{K2015}, for many categories such an analog can be very nontrivial as, for example, for finite graphs, see~\cite{HIK2011}. A more subtle problem concerns the internal direct decompositions. Indeed, in this case the uniqueness of the indecomposable factors up to isomorphism is replaced by the ``absolute'' uniqueness of the set of indecomposable subobjects corresponding to the factors; a good example here is the famous Artin--Wedderburn theorem on semisimple rings. Finally, the third problem is of algorithmic nature, namely, it consists in constructing an efficient algorithm finding a direct product decomposition of a given object into indecomposable subobjects. Two examples here are algorithms for finding the Wedderburn decomposition of associative semisimple algebras~\cite{FR1985} and the direct product decomposition of a finite permutation group~\cite{W2010}. In the present paper, we examine these three problems in relation to coherent configurations.\smallskip

The coherent configurations introduced by D.~Higman in~\cite{H1971} as a tool for studying finite permutation groups, includes finite groups and association schemes  as full subcategories. An analog of the Krull-Schmidt theorem for commutative coherent configurations was found in~\cite{FT1985}. Later, it was generalized to infinite case (with slightly  weaken commutativity condition) in~\cite{Z2005} and to homogeneous (not necessarily commutative) coherent configurations in~\cite{Xu2013}.  We did not find any result on coherent configurations related to the other two problems mentioned above.\smallskip

A coherent configuration on a finite set $\Omega$ can be thought as a finite arc-colored directed complete graph with vertex set $\Omega$, the color classes of which satisfy certain combinatorial conditions (for exact definitions and related material see Section 2). As shown by D.~Higman in~\cite{H1987}, the coherent configurations on $\Omega$ are in one-to-one correspondence with coherent algebras on $\Omega$, which are special semisimple subalgebras of the algebra of all complex matrices whose rows and columns are indexed by the elements of $\Omega$. In the sense of this correspondence, the direct products of coherent configurations are exactly the tensor product of the corresponding coherent algebras. That is why in theory of coherent configurations we use the term “tensor product” rather than “direct product”. \smallskip

A {\it tensor decomposition} of a coherent configuration $\cX$ on~$\Omega$ is a set of  coherent configurations $\cX_i$ on $\Omega_i$, $i=1,\ldots,m$, and  a bijection  $\pi:\Omega\to\Omega_1\times\cdots\times\Omega_m$, which is an isomorphism from $\cX$ to~$\cX_1\otimes\cdots\otimes\cX_m$. Each tensor decomposition of~$\cX$ defines a Cartesian decomposition of~$\Omega$ in the sense of~\cite{BCPS2020,PS2018}. In Section~\ref{110421a}, we introduce the concept of Cartesian decomposition of a coherent configuration and prove Theorems~\ref{250121b} and~\ref{130121c} that establish the following relationship between tensor and Cartesian decompositions of coherent configurations.

\thrml{100421a}
Each tensor decomposition of a coherent configuration $\cX$ defines a uniquely determined Cartesian decomposition of~$\cX$. Conversely, every  Cartesian decomposition of~$\cX$ defines a tensor decomposition of~$\cX$.
\ethrm

The Cartesian decompositions of coherent configurations play the role of ``internal'' tensor decompositions. A maximal Cartesian decomposition of $\cX$ corresponds to a tensor decomposition of~$\cX$  in which every factor is indecomposable with respect to the tensor product. In general, $\cX$ may have several maximal Cartesian decompositions. For example, if the coherent algebra corresponding to $\cX$ is a full matrix algebra, then every Cartesian decomposition of~$\Omega$ is also a Cartesian decomposition of~$\cX$, and the number of distinct maximal Cartesian decompositions of~$\cX$ can be exponential on~$|\Omega|$. In this example, $\cX$ has  irreflexive basis relations of valency~$1$; a coherent configuration is said to be {\it thick} if it has no such basis relations.

\thrml{080619a}
Every thick coherent configuration $\cX$ has exactly one maximal Cartesian decomposition. 
\ethrm

It is not clear whether an analog of the Krull-Schmidt theorem holds for all  coherent configurations. However as a consequence of Theorems~\ref{100421a} and~\ref{080619a}, we have the following statement.

\crllrl{110421c}
The analog of the Krull-Schmidt theorem holds for all  thick coherent configurations.
\ecrllr

Let $G$ be a finite group. The center of the group algebra of $G$ is a homogeneous coherent configuration $\cX$ on $\Omega=G$. The valencies of the basis relations of $\cX$ are exactly the cardinalities of the conjugacy classes of~$G$. Therefore, $\cX$ is thick if and only if the center of~$G$ is trivial. Using Theorem~\ref{100421a}, it is not hard to verify that the decompositions of $G$ into the internal direct product of normal subgroups are in one-to one correspondence with Cartesian decompositions of~$\cX$. Thus by Theorem~\ref{080619a}, we obtain the following statement.

\crllrl{040221a1}
Let $G$ be a finite group with trivial center. Assume that $G=G_1\cdots G_m$ and $G=H_1\cdots H_k$ are two decompositions of $G$ into the internal direct product of  indecomposable  normal subgroups. Then the multisets $\{\{G_1,\ldots,G_m\}\}$ and $\{\{H_1,\ldots,H_k\}\}$ are equal.
\ecrllr

In the last part of the paper, we construct an efficient algorithm finding the maximal  tensor decomposition of a thick coherent configuration~$\cX$. The idea behind the algorithm is somewhat similar to that from~\cite{KN2009}. Using a natural greedy algorithm (Algorithm~A in Section~\ref{270121a}), one can find a nontrivial Cartesian decomposition $P^*$ of the underlying set~$\Omega$. The thickness condition implies that $P^*$ is a refinement of every Cartesian decomposition $P$ of~$\cX$ (Theorem~\ref{170121a}).  Since the number of all possible~$P$ is at most~$|\Omega|$, the maximal  tensor decomposition of $\cX$ can be found by exhaustive search of all~$P$.

\thrml{080619a1}
The maximal  tensor decomposition of a thick coherent configuration of degree~$n$ can be found in polynomial time in~$n$.
\ethrm

The paper is organized as follows. The definitions and basic properties related with coherent configurations  are given in Section~\ref{110421u}. The basic facts concerning Cartesian decompositions of finite sets are discussed in Section~\ref{110421v}. The goal of Section~\ref{110421a} is to introduce Cartesian decomposition of a coherent configuration and then to prove Theorems~\ref{250121b} and~\ref{130121c}. In Sections~\ref{270121a} and~\ref{110421v2}, we prove Theorems~\ref{080619a} and~\ref {080619a1}, respectively.\medskip

The authors are grateful to Professor M. Muzychuk for numerous fruitful discussions on the topic of the paper.

\section{Preliminaries}\label{110421u}

In our presentation of coherent configurations, we mainly follow the monograph~\cite{CP},
where all the details can be found.

\subsection{Notation}\label{120221a}
Throughout the paper, $\Omega$ denotes a finite set. For $\Delta\subseteq \Omega$, the Cartesian product $\Delta\times\Delta$ and its diagonal are denoted by~$\bone_\Delta$ and $1_\Delta$, respectively. For a relation $s\subseteq\bone_\Omega$, we set $s^*=\{(\alpha,\beta): (\beta,\alpha)\in s\}$, $\alpha s=\{\beta\in\Omega:\ (\alpha,\beta)\in s\}$ for all $\alpha\in\Omega$, and define $\grp{s}$ as the minimal equivalence relation on~$\Omega$, containing~$s$. The left and right supports of $s$ are the sets $\Omega_{-}(s)=\{\alpha\in\Omega:\ \alpha s\ne\varnothing\}$ and $\Omega_{+}(s)=\{\alpha\in\Omega:\ \alpha s^*\ne\varnothing\}$, respectively. 
\medskip

The {\it dot product} (sometimes called the composition) of  relations $r$ and $s$ on the same set $\Omega$ is defined as follows:
$$
r\cdot s=\{(\alpha,\beta)\in\bone_\Omega:\alpha r\cap \beta s^*\ne \varnothing\}.
$$
Clearly, $r\cdot s=\varnothing$ if and only if $\Omega_+(r)\cap\Omega_-(s)=\varnothing$, and also $(r\cdot s)^*=s^*\cdot r^*$. We say that  $r$ and $s$ {\it commute} if $r\cdot s=s\cdot r$. \medskip

An equivalence relation is said to be {\it discrete} if each class of it is a singleton.  Given an equivalence relations $e, f$, we set $e\wedge f=e\cap f$ and $e\vee f=\grp{e\cup f}$. The following statement is a part of \cite[Proposition~2.19]{BCPS2020}.

\lmml{180421q}
Let $e$ and $f$ be commuting equivalence relations.  Then $e\cdot f=e\vee f$. Moreover, if $e\wedge f=1_\Omega$, then the intersection of any two classes, one of~$e$ and another one of~$f$, contained in the same class of~$e\cdot f$ is a singleton.
\elmm

Let $e$ be an equivalence relation on $\Omega$ and $\Omega/e$ the set of all classes of~$e$. Then there is  a  surjection
\qtnl{080421a}
\rho_e:\Omega\to\Omega /e,\ \alpha\mapsto\alpha e.
\eqtn 
It induces a natural surjection from the binary relations on $\Omega$ to those on $\Omega/e$, which is also denoted by $\rho_e$.   The following lemma is straightforward.

\lmml{170421d}
Let $e$ be an equivalence relation on $\Omega$ and $\rho=\rho_e$. Then
\nmrt
\tm{1} $\rho(1_\Omega)=1_{\Omega/e}$, $\rho(\bone_\Omega)=\bone_{\Omega/e}$,
\tm{2} if $e$ commutes with $r,s\subseteq\bone_\Omega$, then $\rho(r\cdot s)=\rho(r)\cdot\rho(s)$,
\tm{3} if $e$ commutes with equivalence relations $f,g\subseteq\bone_\Omega$, then $\rho(f\vee g)=\rho(f)\vee\rho(g)$ and $\rho(f\wedge g)=\rho(f)\wedge\rho(g)$.
\enmrt
\elmm

Let $m\ge 1$, and let $\Omega=\Omega_1\times\cdots\times\Omega_m$ be the Cartesian product. The {\it tensor product} of relations $s_1\subseteq\bone_{\Omega_1},\ldots, s_m\subseteq\bone_{\Omega_m}$  is defined to be  
$$
s_1\otimes\cdots\otimes s_m=\{(\alpha,\beta)\in\bone_\Omega:\ (\alpha_i,\beta_i)\in s_i,\ i=1,\ldots,m\},
$$
where $\alpha_i$ and $\beta_i$ are the $i$th components of the $m$-tuples $\alpha$ and $\beta$, respectively. In the matrix language, this means that the adjacency matrix of the tensor product is equal to the Kronecker product of the adjacency matrices of the factors. A relationship between the dot and tensor products is given by the formula
\qtnl{150421a}
(s_1\otimes\cdots\otimes s_m)\cdot (t_1\otimes\cdots\otimes t_m)=(s_1\cdot t_1)\otimes\cdots\otimes (s_m\cdot t_m),
\eqtn
the proof of which is straightforward.\medskip

We define a partial order $\preceq$ on the sets of relations on the same set by setting  $S\preceq T$ if every relation of~$T$ is contained in some relation of~$S$.

\subsection{Coherent configurations}

Let $\Omega$ be a finite set and $S$ a partition of $\Omega^2$; in particular, the elements of $S$ are treated as binary relations on~$\Omega$. A pair $\mathcal{X}=(\Omega,S)$ is called a \emph{coherent configuration} on $\Omega$ if the following conditions are satisfied:
\nmrt
\tm{C1}  $1_\Omega$ is the union of some relations of~$S$,
\tm{C2} $s^*\in S$ for all $s\in S$, 
\tm{C3} given $r,s,t\in S$ the number $c_{rs}^t=|\alpha r\cap \beta s^{*}|$ does not depend on $(\alpha,\beta)\in t$. 
\enmrt
The number $|\Omega|$ is called the {\it degree} of~$\cX$. We say that $\cX$ is {\it trivial} if  $S$ consists of~$1_\Omega$ and its complement (unless $\Omega$ is not a singleton), {\it homogeneous} if $1_\Omega\in S$, and {\it commutative} if $c_{rs}^t=c_{sr}^t$ for all $r,s,t$. \medskip 

An {\it isomorphism} from $\cX$ to a coherent configuration  $\cX'=(\Omega', S')$ is a bijection $\pi: \Omega\rightarrow \Omega'$ such that for each $s\in S$, the relation $\pi(s)=\{(\alpha^\pi, \beta^\pi): (\alpha, \beta)\in s\}$ belongs to~$S'$. In this case,  we say that $\cX$ and $\cX'$ are {\it isomorphic} and write $\cX\cong\cX'$.\medskip

Let $G$ be a permutation group on~$\Omega$. Denote by $S$ the set of all orbits $(\alpha,\beta)^G$ of the induced action of~$G$ on~$\Omega\times\Omega$, $\alpha,\beta\in \Omega$. Then the pair $(\Omega,S)$ is a coherent configuration; we say that it is the coherent configuration associated with $G$. 

\subsection{Relations}  The elements of $S=S(\cX)$ are called {\it basis relations} of the coherent configuration~$\cX$. We say that $s\in S$ is {\it reflexive} if $s\subseteq 1_\Omega$, and {\it irreflexive} otherwise. The unique basis relation containing the pair~$(\alpha,\beta)$ is denoted by~$r(\alpha,\beta)$. Thus, $r(\alpha,\beta)=r(\beta,\alpha)^*$. \medskip

A {\it fiber} of $\cX$ is a set $\Delta\subseteq\Omega$ such that $1_\Delta\in S$. By condition~(C1), the set of all fibers form a partition of~$\Omega$. Moreover,  the sets $\Omega_-(s)$ and $\Omega_+(s)$ are fibers of $\cX$ for all $s\in S$. In particular, $S$ is partitioned into the disjoint union of the sets
$$
S_{\Delta,\Gamma}=\{s\in S:\ \Omega_-(s)=\Delta,\ \Omega_+(s)=\Gamma\},
$$
where $\Delta$ and~$\Gamma$ run over the fibers of~$\cX$.\medskip

Any union of basis relations is called a {\it relation} of $\cX$.  The set of all of them is closed with respect to intersections, unions, and the dot product. Together with condition~(C1) this implies that if $s$ is a nonempty relation of $\cX$, then so are $1_{\Omega_-(s)}$ and~$1_{\Omega_+(s)}$.  Moreover, if $s\in S$, then the latter relations are basis.\medskip

Let $s\in S$ and $t=1_{\Omega_-(s)}$. The integer $n_s=c_{ss^*}^t$  is positive and is equal to $|\alpha s|$ for all $\alpha\in \Omega_-(s)$; it  is called the {\it valency} of $s$.  We say that $s$ is {\it thin} if $n_{s^{}}=n_{s^*}=1$; in particular, every reflexive basis relation is thin. An important property of a thin relation $s\in S$ is that given $r\in S$ each of the dot products $r\cdot s$ and $s\cdot r$ is either a basis relation  or empty. A coherent configuration is said to be {\it thick} if it contains no  irreflexive thin relations.

\subsection{Parabolics and quotients} An equivalence relation which is a relation of the coherent configuration~$\cX$ is called a {\it parabolic} of $\cX$.  The set of all parabolics of~$\cX$ is denoted by $E=E(\cX)$. It is a lattice the minimal and maximal elements of which are {\it trivial} parabolics $1_\Omega$ and $\bone_\Omega$; the meet   $\wedge$ and join~$\vee$ are exactly the operations  defined in Subsection~\ref{120221a}. Note that if parabolics~$e$ and $f$ commute, then $e\vee f=e\cdot f=f\cdot e$.\medskip

We define a relation $\perp$ on the set $E$ as follows: $e\perp f$ if for any two basis relations $x\subseteq e$ and $y\subseteq f$ we have 
\qtnl{060221c}
x\cdot y\ne\varnothing\ \, \Rightarrow\ \,x\cdot y\in S.
\eqtn
This relation is not necessarily reflexive, but is symmetric. Indeed, if $y\cdot x\ne\varnothing$, then $x^*\cdot y^*\ne\varnothing$ and $y\cdot x=(x^*\cdot y^*)^*\in S$. From the definition, it immediately follows that given parabolics $e', f'\in E$, we have
\qtnl{070221a}
e'\subseteq e,\ f'\subseteq f,\ e\perp f\qquad\Rightarrow\qquad e'\perp f'.
\eqtn 

\lmml{230421i}
If $\cX$ is thick and $e,f\in E$ are such that $e\perp f$, then $e\wedge f=1_\Omega$.
\elmm
\proof  Assume on the contrary  that $e\wedge f\ne 1_\Omega$. Then  there is an irreflexive basis relation $s\subseteq e\wedge f\subseteq e$.  Since $s^*\subseteq e\wedge f\subseteq f$,  $s\cdot s^*\ne\varnothing$, and $s^*\cdot s\ne\varnothing$, the condition $e\perp f$ implies that $s\cdot s^*,s^*\cdot s\in S$. Since $1_{\Omega_-(s)}\subseteq s\cdot s^*$ and $1_{\Omega_+(s)}\subseteq s^*\cdot s$, we have 
$$
s\cdot s^*=1_{\Omega_-(s)}\qaq s^*\cdot s=1_{\Omega_+(s)}.
$$
Consequently, $n_{s^{}}=n_{s^*}=1$, which is impossible, because $\cX$ is thick.\eprf\medskip

We say that parabolics $e,f\in E$ {\it strongly commute} if for any two basis relations $x\subseteq e$ and $y\subseteq f$, we have 
$$
x\cdot f=f\cdot x\qaq e\cdot y=y\cdot e.
$$
Again the relation ``to be strongly commute'' is not reflexive but is symmetric. Clearly, if $e$ and $f$ strongly commute, then $e$ and $f$ commute. But the reverse statement is not true.\medskip

Let $e$ be a parabolic of~$\cX$. Denote by $S_{\Omega/e}$ the set of all relations $\rho_e(s)$, $s\in S$. Then the pair
$$
\cX_{\Omega/e}=(\Omega/e,S_{\Omega/e})
$$
is a coherent configuration; it is called the {\it quotient} of~$\cX$ modulo~$e$. The mapping~$\rho_e$ induces a surjection from the relations (respectively, parabolics) of $\cX$ to the relations (respectively, parabolics) of $\cX_{\Omega/e}$.

\subsection{Tensor product}
Let $\cX_i=(\Omega_i, S_i)$ be a coherent configuration, $i=1,\ldots,m$. Denote by $S_1\otimes\cdots\otimes S_m$ the set of all tensor products  $s_1\otimes\cdots\otimes s_m$,  where $s_i\in S_i$ for all~$i$. The pair 
$$
\cX_1\otimes\cdots\otimes \cX_m:=(\Omega_1\times\cdots\times \Omega_m, S_1\otimes\cdots\otimes S_m)
$$
is a coherent configuration; it is called  the {\it tensor product} of $\cX_1,\ldots,\cX_m$.  We say that the tensor product is nontrivial if each factor is of degree at least~$2$, and trivial otherwise. A coherent configuration is said to be {\it indecomposable} if it is not isomorphic to a nontrivial tensor product of  coherent configurations.

\section{Cartesian decompositions}\label{110421v}

\subsection{Atomic Cartesian decomposition of a set.}
In what follows, $m$ is a positive integer and $P=\{e_1,\ldots,e_m\}$ is a set of equivalence relations on~$\Omega$. The index set $\{1,\ldots,m\}$ is denoted by $M=M(P)$.\medskip

Assume that the relations $e_i\in P$ commute  pairwise. For every $I\subseteq M$, the dot product~$P_I$ of all $e_i$, $i\in I$, is again an equivalence relation  (Lemma~\ref{180421q}); when~$I$ is empty, $P_I=1_\Omega$.  Moreover, if  $I,J\subseteq M$, then
\qtnl{010221b1}
P_I\vee P_J=P_{I\cup J}.
\eqtn
The $P$-complement of $e=P_I$  is defined to be $e'=P_{M\setminus I}$; so, $e\vee e'=P_M$. Thus, 
$L(P)=\{P_I:\ I\subseteq M\}$ is a (join) semilattice in which every element has a complement, and the minimal and maximal elements equal $P_{\varnothing}$ and ~$P_M$, respectively. Note that, in general, $L(P)$ is not a lattice, see \cite[Example 2.21]{BCPS2020}.\smallskip

 The following statement is probably well-known, but the authors are not aware of the publication where this fact was originally established.

\thrml{020421a}
Let $P=\{e_i:\ i\in M\}$ be a  set of pairwise commuting equivalence relations on $\Omega$. Assume that
$e^{}_i\wedge e_i'=1_\Omega$ for all $i\in M$. Then  for all $I,J\subseteq M$, 
$$
P_I\wedge P_J=P_{I\cap J}.
$$
In particular, $L(P)$ is a lattice isomorphic to the Boolean lattice $2^M$.
\ethrm
\proof  First, let $I\cap J=\varnothing$. In this case,  $P_J\subseteq P_{M\setminus I}$, and  we may assume that $J=M\setminus I$, and also $|I|\le |J|$. The proof goes by induction on $|I|$. When $|I|=1$, the statement follows from the theorem assumption.\medskip

Let $|I|\ge 2$. Assume on the contrary that the equivalence relation $f:=P_I\cap P_{M\setminus I}$ is not discrete. Then $P_I$ is not discrete; let $i\in I$ be such that $e=e_i$ is not discrete. Put  $\ov e=P_{I\setminus \{i\}}$ and $g=P_{M\setminus I}$;  in particular, $f=(e\vee \ov e)\wedge g$. Clearly, $f\subseteq e\vee\ov e$ and $f\subseteq g$. We claim that
\begin{figure}[t]
$\xymatrix@R=8pt@C=15pt@M=0pt@L=2pt{
\VRTB{e\vee \ov e}  & & \VRTB{g\vee \ov e} & & \\
& & & & \\
& & \VRTB{f\vee \ov e}\ar@{--}[lluu]\ar@{-}[uu]     & & \\
& & & & \\
\VRTB{e} \ar@{-}[uuuu]  & &\VRTB{\ov e}\ar@{-}[uu] & & \VRTB{g} \ar@{-}[lluuuu] \\
& & & \VRTB{f}\ar@{--}{luu}\ar@{-}[ru]\ar@{--}[luuu]&\\
& & \VRTB{1}\ar@{-}[lluu] \ar@{-}[uu] \ar@{-}[ru]& &\\
}$
\caption{The diagram of the partial order from Theorem~\ref{020421a}.}\label{fig2}
\end{figure}
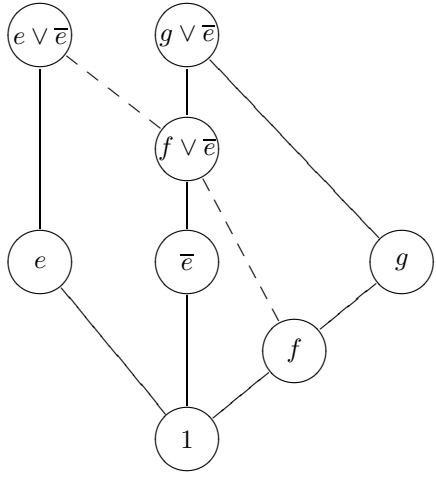
\qtnl{170421a}
\ov e\,\subsetneq\, f\vee \ov e\,\subsetneq\, e\vee \ov e\qaq (f\vee \ov e)\wedge e=1_\Omega.
\eqtn
Indeed, if $\ov e=f\vee \ov e$, then $f\subseteq \ov e$, and hence $f\subseteq \ov e\wedge g$. However, $\ov e=P_{I\setminus \{i\}}$ and then by induction, $\ov e\wedge g=1_\Omega$ implying $f=1_\Omega$, a contradiction. This proves that $\ov e\,\subsetneq\, f\vee \ov e$.  Next, by the theorem assumption, we have $e\wedge e'=1_\Omega$ and hence
$$
(f\vee \ov e)\wedge e\subseteq (g\vee \ov e)\wedge e=  e'\wedge e=1_\Omega.
$$ 
This proves  the second formula in~\eqref{170421a}. Finally, if $f\vee \ov e= e\vee \ov e$, then $(f\vee \ov e)\wedge e=( e\vee \ov e)\wedge e=e\ne 1_\Omega$ in contrast to what we just proved.  Thus the diagram of the partial order between the defined equivalence relations is as in Fig.~\ref{fig2}\medskip

The first formula in~\eqref{170421a} implies that there are two distinct classes $\Delta_1$, $\Delta_2$ of $\ov e$, contained in the same class $\Delta$ of  $f\vee\ov e$. Since $\ov e$ and $e$ commute and $\Delta$ is contained in a class of~$e\vee\ov e$, there are points $\delta_1\in \Delta_1$ and $\delta_2\in\Delta_2$ such that $(\delta_1,\delta_2)\in e$ (Lemma~\ref{180421q}). Since $\delta_1,\delta_2\in\Delta_1\cup\Delta_2\subseteq\Delta$, this shows that $(\delta_1,\delta_2)\in f\vee \ov e$. Thus, $ (f\vee \ov e)\wedge e\ne 1_\Omega$, in contrast to the second formula in~\eqref{170421a}. This completes the proof for the case $I\cap J=\varnothing$.\medskip

Now we turn to the case $I\cap J\ne\varnothing$. Without loss of generality, we may assume that $|M|>1$.  Let $i\in I\cap J$ and put 
$$
\rho=\rho_{e_i}\qaq \ov P=\{\rho(e_k):\ k\in \ov M\},
$$ 
where for all $K\subseteq M$,  we set $\ov K=K\setminus\{i\}$.  By Lemma~\ref{170421d}(2), the relations of  $\ov P$, pairwise commute. Furthermore, given  $\ov K\subseteq \ov M$, we have
\qtnl{170421w}
\rho(P_K)=\ov P_{\ov K}.
\eqtn
Indeed, the left-hand side obviously contains the right-hand side. On the other hand, if this inclusion is strict, then there exists $l\not\in K\cup\{i\}$ such that $\rho(e_l)\subseteq \rho(P_K)$. Then $e_l\subseteq e_l\cdot e_i\subseteq P_{K\cup\{i\}}\subseteq e'_l$ and hence $e_l\subseteq e^{}_l\wedge e'_l$, in contrast to the  theorem assumption.\medskip

Thus, the set $\ov P$ satisfies the theorem assumption. Using Lemma~\ref{170421d}(3), formula~\eqref{170421w}, and  induction on $|P|$,  we obtain
$$
\rho(P_I\wedge P_J)=\rho(P_I)\wedge\rho(P_J)=\ov P_{\ov I}\cap\ov P_{\ov J}=\ov P_{\ov I\cap\ov J}.
$$ 
Since $i\in I\cap J$, the $\rho$-preimage of the left-hand side equals $P_I\cap P_J$, and $\ov I\cap \ov J=\ov{ I\cap J}$. By formula~\eqref{170421w}, this implies that
$$
P_I\wedge P_J=\rho^{-1}(\rho(P_I\wedge P_J))=\rho^{-1}(\ov P_{\ov I\cap\ov J})=\rho^{-1}(\rho(P_{I\cap J})=P_{I\cap J},
$$
as required.\eprf

\crllrl{040421a}
In the condition of Theorem~\ref{020421a}, assume that $P_M=\bone_\Omega$. Let $I\subseteq M$,  and let~$\Delta'_i$ be a class of $e'_i$ for each $i\in I$. Then the intersection of $\Delta_i'$, $i\in I$, is a class of~$P_{M\setminus I}$.
\ecrllr
\proof For $|I|=1$, the statement is trivial. Let $|I|>1$, $i\in I$, and $J=I\setminus\{i\}$. By induction, the intersection $\Gamma$ of all classes $\Delta'_j$, $j\in J$, is a class of $P_{M\setminus J}$. We need to verify that  $\Delta'_i\cap\Gamma$ is a class of $P_{M\setminus I}$. However, either $\Delta'_i\cap\Gamma=\varnothing$  or $\Delta'_i\cap\Gamma$ is a class of the equivalence relation
$$
e'_i\wedge P_{M\setminus J}=P_{M\setminus \{i\}}\wedge P_{M\setminus J}=P_{M\setminus I},
$$
where the last equality follows from Theorem~\ref{020421a}. It remains to note that the former case is impossible. Indeed,  the equivalence relations $e'_i$ and~$P_{M\setminus J}$ commute and $e'_i\vee P_{M\setminus J}\supseteq P_{M\setminus\{i\}}\vee P_{\{i\}}=P_M=\bone_\Omega$. Thus each class of $e'_i$ intersects each class of~$P_{M\setminus J}$ (Lemma~\ref{180421q}).\eprf\medskip.

From  Corollary~\ref{040421a} for $I=M$, it follows that if $P_M=\bone_\Omega$ and $|M|>1$, then $\{e'_i:\ i\in M\}$  is a {\it Cartesian decomposition} of $\Omega$ in the sense of~\cite[Section~3.1]{BCPS2020}.  In particular, there is a well defined bijection 
\qtnl{050421f}
\pi_P:\Omega/e'_1\times\cdots\times\Omega/e'_m\to \Omega,\ (\Delta'_1,\ldots,\Delta'_m)\mapsto\delta,
\eqtn 
where $\delta$ is a unique point of the intersection $\Delta'_1\cap\cdots\cap\Delta'_m$ (Corollary~\ref{040421a}). It immediately follows that 
\qtnl{190421a}
\rho_{e'_i}(\pi_P(\Delta'_1,\ldots,\Delta'_m))=\Delta'_i,\quad i\in M,
\eqtn
and
\qtnl{080221c}
|\Omega|= \prod_{i\in M} |\Omega/e'_i|.
\eqtn
Moreover, Theorem~\ref{020421a} shows that $L(P)$ is the Cartesian lattice defined by this decomposition (again in the sense of~\cite[Section~3.1]{BCPS2020}), and $P$ consists of the atoms of this lattice.

\dfntnl{080221a}
A nonempty set $P$ of nondiscrete equivalence relations on $\Omega$ is called an atomic Cartesian decomposition of $\Omega$  if they pairwise commute and for every $e\in P$,
\qtnl{180421v}
e\wedge e'=1_\Omega\qaq e\vee e'=\bone_\Omega. 
\eqtn
The atomic Cartesian decomposition is {\it trivial} if $|P|=1$.
\edfntn

The mapping $P\mapsto P'=\{e':\ e\in P\}$ is a poset anti-isomorphism between the atomic and ordinary Cartesian decompositions of the same set $\Omega$. An advantage to use atomic Cartesian decomposition comes from the fact that condition~\eqref{180421v} is sometimes more easier to check (at least from an algorithmic  point of view) than the conclusion of Corollary~\ref{040421a} for $I=M$. 

\subsection{Standard atomic Cartesian decomposition.}

Let  $\Omega=\Omega_1\times\cdots\times\Omega_m$, where $|\Omega_i|>1$ for all $i\in M$.  Denote by~$e_i$ the equivalence relation on $\Omega$ defined by the equalities at all coordinates except for the $i$th one, i.e.,
$$
e_i=1_{\Omega_1}\otimes\cdots\otimes 1_{\Omega_{i-1}} \otimes \bone_{\Omega_i}\otimes 1_{\Omega_{i+1}}\otimes\cdots \otimes 1_{\Omega_m},\quad i\in M.
$$
These equivalence relations are obviously pairwise commuting. Put $P=\{e_i:\ i\in M\}$. Then the $P$-complement of $e_i$  is given by formula
$$
e'_i=\bone_{\Omega_1}\otimes\cdots\otimes \bone_{\Omega_{i-1}} \otimes 1_{\Omega_i}\otimes \bone_{\Omega_{i+1}}\otimes\cdots \otimes \bone_{\Omega_m},\quad i\in M.
$$
It easily follows that $P$ is an atomic Cartesian decomposition of~$\Omega$; it is said to be a {\it standard} one. 

\lmml{060421a}
Let $P$ be an arbitrary atomic Cartesian decomposition and $\pi=\pi_P$ the bijection~\eqref{050421f}. Then $\pi^{-1}(P)$ is the standard atomic Cartesian decomposition of the~$\dom(\pi)$.
\elmm
\proof  The set $\dom(\pi)$ consists of the tuples $\Delta'=(\Delta'_1,\ldots,\Delta'_m)$, where $\Delta'_i$ is a class of $e'_i$, $i\in M$. The definitions of $e_i$ and $\pi$ imply that 
$$
\pi^{-1}(e_i)=\{(\Delta',\Gamma'):\ \Delta'_j=\Gamma'_j\ \text{ for all } j\in M\setminus\{ i\}\}
$$
is exactly the equivalence relation on $\dom(\pi)$ defined by the equalities at all coordinates except for the $i$th one.\eprf\medskip

Let  $I=\{i_1,\ldots,i_k\}$, where $1\le i_1<i_2<\cdots<i_k\le m$. The projection of $\Omega$ with respect to $I$ is defined to be the surjection
$$
\pr_I:\Omega\to \Omega_{i_1}\times\cdots\times\Omega_{i_k},\ (\alpha_1,\ldots,\alpha_m)\mapsto (\alpha_{i_1},\ldots,\alpha_{i_k}).
$$
It is easily seen that $\pr_I(P):=\{\pr_I(e_{i_1}),\ldots,\pr_I(e_{i_k})\}$  is the standard atomic Cartesian decomposition of $\pr_I(\Omega)=\Omega_{i_1}\times\cdots\times\Omega_{i_k}$.

\subsection{Properties of atomic Cartesian decompositions.}
The atomic Cartesian decompositions on $\Omega$ are partially ordered with respect to the relation~$\preceq$. The following lemma shows that nontrivial upper ideals of this partially ordered set admit a good description.

\lmml{040421c}
Let $P$ be an atomic Cartesian decomposition of $\Omega$. Then
\nmrt
\tm{1} for every  partition $\cR$ of~$M$, the set $P^\cR=\{P_I:\ I\in\cR\}$ is  an atomic Cartesian decomposition of $\Omega$,
\tm{2} if $R\preceq P$ is an atomic Cartesian decomposition of $\Omega$, then $R=P^\cR$ for some
partition $\cR$ of $M$.
\enmrt
\elmm
\proof Statement (1) immediately follows from formula~\eqref{010221b1} and Theorem~\ref{020421a}. To prove statement~(2), let $f\in R$ and 
$$
M(f)=\{i\in M:\ e_i\subseteq f\}.
$$  
Since $R\preceq P$, every element of $P$ is contained in either~$f$, or in some element of $R$, contained in the $R$-complement~$f'$ of~$f$. Thus,
\qtnl{040421q}
P_{M(f)}\subseteq f\qaq P_{M\setminus M(f)}\subseteq f'.
\eqtn
In particular, the nonempty sets $M(f)$, $f\in R$, are pairwise disjoint and  form a partition~$\cR$ of~$M$. It suffices to verify that $P_{M(f)}=f$ for all $f\in R$.\medskip

By statement (1) applied to $P$ and $\cR$,  each of the sets $\{P_{M(f)},P_{M\setminus M(f)}\}$ and $\{f,f'\}$ is an atomic Cartesian decomposition of $\Omega$. By formula~\eqref{080221c}, this yields
$$
|\Omega/P_{M(f)}|\,\cdot\,|\Omega/P_{M\setminus M(f)}|=|\Omega|=|\Omega/f|\,\cdot\,|\Omega/f'|.
$$
On the other hand, by formula~\eqref{040421q}, we have
$$
|\Omega/P_{M(f)}|\ge |\Omega/f|\qaq |\Omega/P_{M\setminus M(f)}\ge |\Omega/f'|.
$$
These two formulas together imply that 
$$
|\Omega/P_{M(f)}|=|\Omega/f|\qaq|\Omega/P_{M\setminus M(f)}|=|\Omega/f'|,
$$ 
which, in particular, means that $P_{M(f)}=f$, as required.\eprf\medskip

Let $P$ and $Q$ be atomic  Cartesian decompositions of~$\Omega$.  In general, the set
$$
P\cap Q=\{e\wedge f:\  e\in P,\ f\in Q,\ e\wedge f\ne 1_\Omega\}
$$
is not an atomic Cartesian decomposition 
The following statement shows that this never happens if there is an atomic  Cartesian decomposition of~ $\Omega$, which is a common refinement of~$P$ and~$Q$.

\lmml{260121a}
Let $Q$ and $R$ be  atomic Cartesian decompositions of $\Omega$. Assume that there is an atomic  Cartesian decomposition~$P$ of~ $\Omega$ such that $ Q\cap R\preceq P$. Then $Q\cap R$ is also an atomic Cartesian decomposition of $\Omega$.
\elmm
\proof   We have $Q\preceq Q\cap R\preceq P$ and $R\preceq Q\cap R\preceq P$. By Lemma~\ref{040421c}(2),  this implies that there exist  partitions $\cQ,\cR$ of $M$ such that
$$
Q=P^\cQ\qaq R=P^\cR.
$$ 
For any $f\in Q$ and any $g\in R$, there are $I\in\cQ$ and $J\in\cR$ for which $f=P_I$ and $g=P_J$.  By Theorem~\ref{020421a}, we have 
$$
f\wedge g=P_I\wedge P_J=P_{I\cap J}.
$$ 
The nonempty sets $I\cap J$ form a partition $\cP$ of the set $M$. Consequently,  $Q\cap R=P^\cP$ and we are done by Lemma~\ref{040421c}(1).  \eprf

\section{Cartesian decompositions and tensor products}\label{110421a}

Unless otherwise stated, in what follows we always assume that $\cX=(\Omega,S)$ is a coherent configuration and $E=E(\cX)$.\medskip

\dfntnl{229421f}
An (atomic) Cartesian decomposition $P$ of $\Omega$ is called the (atomic) Cartesian decomposition of $\cX$, if $P\subseteq E$ and each $e\in P$ strongly commutes with $e'$, and also $e\perp e'$.
\edfntn

Clearly, the poset anti-isomorphism between the atomic and ordinary Cartesian decompositions of a set~$\Omega$ induces a poset anti-isomorphism between the atomic and ordinary Cartesian decompositions of $\cX$. Thus in the sequel we deal with atomic Cartesian decompositions only.\medskip

Let $P$ be an atomic Cartesian decomposition of $\cX$. Clearly, the set $\{e,e'\}$ is also an atomic Cartesian decomposition of $\cX$ for all $e\in P$.  Furthermore, for any $I\subseteq M$, the parabolic $P_{M\setminus I}$
commutes with all relations of~$P$. By Lemma~\ref{170421d}, this implies that $\rho_e(P)=\{\rho_e(f):\ f\in P\}$ is an atomic Cartesian decomposition of the coherent configuration~$\cX_{\Omega/e}$, where $e=P_{M\setminus I}$.

\lmml{180421a}
Let $P$ be an atomic Cartesian decomposition of $\cX$, and let $\Delta$, $\Gamma$ be fibers of~$\cX$. Then
\nmrt
\tm{1}  $\Delta=\Gamma$ if and only if  $\rho_{e'}(\Delta)=\rho_{e'}(\Gamma)$ for all  $e\in P$,
\tm{2}  the restriction of $\rho_{e'}$ to $\{s\in S_{\Delta,\Gamma}:\ s\subseteq e\}$ is injective for  all $e\in P$.
\enmrt
\elmm
\proof  (1) Let $\pi=\pi_P$. Recall that $\{e,e'\}$ is a Cartesian decomposition of~$\Omega$. We claim that 
$$
\rho_{e^{}}(\Delta)\times\rho_{e'}(\Delta)=\pi^{-1}(\Delta).
$$ 
Indeed, the inclusion $\supseteq$ is obvious. To prove the reverse inclusion, let $\Lambda\in\rho_e(\Delta)$ and $\Lambda'\in\rho_{e'}(\Delta)$, i.e., there are points $\alpha,\alpha'\in\Delta$ such that $\rho_e(\alpha)\in\Lambda$ and $\rho_{e'}(\alpha')\in\Lambda'$. By the definition of $\pi$, $\pi(\Lambda,\Lambda')$ is equal to the unique point $\beta$ of the set~$\Lambda\cap\Lambda'$. We need to verify that $\beta\in\Delta$.\medskip

Note that $r(\alpha,\beta)\subseteq e$, because both $\alpha$ and $\beta$ belong to the class $\Lambda$ of $e$, and similarly, $r(\beta,\alpha')\subseteq e'$. Since $e\perp e'$, we have 
$$
r(\alpha,\alpha')=r(\alpha,\beta)\cdot r(\beta,\alpha').
$$
Furthermore,  $e$ and $e'$ strongly commute. Consequently, there exists a basis relation $s'\subseteq e'$ such that
$$
r(\alpha,\alpha')=s'\cdot r(\alpha,\beta).
$$
It follows that $\beta\in\Omega_+(s)$, where $s=r(\alpha,\alpha')$. It remains to note that $\Omega_+(s)=\Delta$, because $\alpha'\in\Delta$. The claim is proved.\medskip

Now from the claim, it easily follows that for every fiber $\Lambda$ of $\cX$, we have 
$$
\pi(\Lambda)=\rho_{e'_1}(\Lambda)\times\cdots\times\rho_{e'_m}(\Lambda).
$$ 
Therefore if $\pi_{e'_i}(\Delta)=\pi_{e'_i}(\Gamma)$ for all $i\in M$, then $\pi(\Delta)=\pi(\Gamma)$. Since $\pi$ is a bijection, we obtain $\Delta=\Gamma$.\medskip

(2)  Let $s,t\in S_{\Delta,\Gamma}$ be such that $\rho_{e'}(s)=\rho_{e'}(t)$.  Then there exist classes $\Lambda',\Lambda''$ of $e'$ and  points 
\qtnl{230421q}
\alpha_s,\alpha_t\in\Lambda'\cap\Delta,\quad \beta_s,\beta_t\in\Lambda''\cap\Gamma,
\eqtn
for which $(\alpha_s,\beta_s)\in s$ and $(\alpha_t,\beta_t)\in t$. Assume that $r(\alpha_s,\beta_s),r(\alpha_t,\beta_t,)\subseteq e$. Since $e\perp e'$ and $e,e'$  strongly commute, we have
$$
r(\alpha_t,\alpha_s)\cdot r(\alpha_s,\beta_s)=r(\alpha_s,\beta_s)\cdot s'=s\cdot s'.
$$
for some $s'\in S$ contained in $e'$. It follows that the set $\alpha_t r(\alpha_s,\beta_s)=\alpha_t s$ intersects the class $\Lambda''$. Thus, we may assume that $\alpha_s=\alpha_t:=\alpha$. Now,
$$
(\alpha,\beta_s)\in s\subseteq e \qaq (\alpha,\beta_t)\in t\subseteq e. 
$$
Therefore, $(\beta_s,\beta_t)\in e$. On the other hand, $(\beta_s,\beta_t)\in e'$ by formula~\eqref{230421q}. Since $e\wedge e'=1_\Omega$, this implies that $\beta_s=\beta_t:=\beta$. Thus, $s$ and $t$ contain a common pair $(\alpha,\beta)$ and hence $s=t$. \eprf\medskip

The  following two theorems establish a one-to-one correspondence between the tensor decompositions  and atomic Cartesian decompositions of~$\cX$.

\thrml{250121b}
Let $\cX=\cX_1\otimes\cdots\otimes\cX_m$, where $\cX_i$ is a coherent configuration on~$\Omega_i$, $i=1,\ldots,m$. Then the standard atomic Cartesian decomposition~$P=\{e_1,\ldots,e_m\}$ of $\Omega=\Omega_1\times\cdots\times\Omega_m$  is an atomic Cartesian decomposition of $\cX$. Moreover, 
$$
\cX_i\cong\cX_{\Omega/e'_i},\quad i=1,\ldots,m.
$$  
\ethrm
\proof  Clearly,  $P\subseteq E$. Next,  $\cX=\cX^{}_i\otimes\cX'_i$ for all~$i$,  where 
$$
\cX'_i=\cX_1\otimes\cdots\otimes\cX_{i-1}\otimes\cX_{i+1}\otimes \cdots\otimes\cX_m
$$
is a coherent configuration on
$$
\Omega'_i=\Omega_1\times\cdots\times\Omega_{i-1}\times\Omega_{i+1}\times \cdots\times\Omega_m.
$$  
Let $x\subseteq e_i$ and $x'\subseteq e'_i$ be arbitrary basis relations of~$\cX$. There exist sets $\Delta^{}_i\subseteq \Omega^{}_i$  and $\Delta'_i\subseteq\Omega'_i$ (which are, in fact, fibers of $\cX^{}_i$ and~$\cX'_i$, respectively), and relations $s^{}_i\in S(\cX_i)$ and $s'_i\in S(\cX'_i)$ such that
$$
x=s^{}_i\otimes 1_{\Delta'_i}\qaq x'=1_{\Delta^{}_i}\otimes s'_i.
$$
Assume that $x\cdot x'\ne\varnothing$. Then $\Delta_i=\Omega_+(s_i)$ and $\Delta'_i=\Omega_-(s'_i)$. By formula~\eqref{150421a}, we obtain
$$
x\cdot x'=(s^{}_i\otimes 1_{\Delta'_i})\cdot(1_{\Delta^{}_i}\otimes s'_i)=(s^{}_i\cdot 1_{\Delta_i})\otimes(1_{\Delta'_i}\cdot s'_i)=s^{}_i\otimes s'_i.
$$
Now  the relation on the right-hand side belongs to $S$. Thus, $e_i\perp e'_i$. Next,  $x''=1_{\Omega_-(s_i)}\otimes s'_i$ is a basis relation of $\cX$, contained in~$e'_i$. Moreover, again by formula~\eqref{150421a}, we have
$$
x\cdot x'=s^{}_i\otimes s'_i=(1_{\Omega_-(s_i)}\otimes s'_i)\cdot (s^{}_i\otimes 1_{\Delta'_i})=x''\cdot x.
$$
Hence, $x\cdot e_i'=e_i'\cdot x$. Similarly, one can verify that $e_i\cdot y=y\cdot e_i$. Therefore, $e^{}_i$ and~$e'_i$ strongly commute. Thus, $P$  is an atomic Cartesian decomposition of $\cX$.  The  statement about the isomorphism follows from~\cite[Theorem 3.2.9]{CP} applied to the tensor product $\cX=\cX^{}_i\otimes\cX'_i$.\eprf\medskip 

The reverse statement to  Theorem~\ref{250121b},   is given below.

\thrml{130121c}
Let  $\cX$ be a coherent configuration on $\Omega$, $P$ an atomic  Cartesian decomposition of $\cX$, $\pi_P$ the bijection \eqref{050421f}, and 
$$
\cX'=\bigotimes\limits_{e\in P}\cX_{\Omega/e'}.
$$
Then $\pi_P$ is an isomorphism from $\cX'$ to $\cX$, and  $P=\pi_P(P')$, where $P'$ is the standard  atomic  Cartesian decomposition of $\cX'$.
\ethrm
\proof  Let $P=\{e_1, \ldots, e_m\}$ and $\pi=\pi_P$. By Lemma~\ref{060421a}, we have $P=\pi(P')$; in particular, $P'=\{f_1,\ldots,f_m\}$, where $f_i=\pi(e_i)$, $i\in M$. To prove that $\pi$  is an isomorphism from $\cX'$ to $\cX$, we verify that $\pi(s')\in S$ for all basis relations $s'$ of~$\cX'$. Let 
$$
s'=s'_1\otimes\cdots\otimes s'_m,
$$ 
where $s'_i\in S(\cX_i)$, $i\in M$.\medskip 

First, assume that $s'$ is reflexive. In this case, $\pi(s')\subseteq 1_\Omega$, and $s'_i=1_{\Delta'_i}$ for all $i\in M$, where $\Delta'_i$ is a fiber of the coherent configuration $\cX_{\Omega/e'_i}$. By formula~\eqref{190421a},
$$
\rho_{e'_i}(\pi(\Delta'_1,\ldots,\Delta'_m))=\Delta'_i,\quad i\in M.
$$
It follows that $\pi(\Delta'_1,\ldots,\Delta'_m)$ is equal to the intersection $\Delta$ of the $\rho_{e'_i}$-preimages of the sets~$\Delta'_i$, $i\in M$. Each of these preimages is a union of some fibers of~$\cX$. Moreover, $\Delta$ cannot contain two distinct fibers of~$\cX$ by Lemma~\ref{180421a}(1).  Thus, $\Delta$ is a fiber of~$\cX$ and $\pi(s')= 1_\Delta\in S$. In fact, we proved that $\pi$ takes the fibers of~$\cX'$ to the fibers of~$\cX$.\medskip

Now let $s'$ be irreflexive. We use induction on the number $k$ of all $i\in M$ for which the relation $s'_i$ is irreflexive. Let $k=1$. Then $s'\subseteq f_i$ for exactly one~$i\in M$. Therefore, $\pi(s')\subseteq \pi(f_i)=e_i$. Since also $\pi(s')\subseteq \rho_{e'}^{-1}(s'_i)$, we obtain
$$
\pi(s')\subseteq e_i\cap  \rho_{e'_i}^{-1}(s'_i).
$$
The reverse inclusion follows from the definition of~$\pi$. Thus, $\pi(s')$ equals the right-hand side of the above formula. On the other hand, $s'\in S'_{\Delta',\Gamma'}$, where $S'=S(\cX')$ and $\Delta'$, $\Gamma'$ are fibers of~$\cX'$. From the above paragraph, it follows that $\Delta=\pi(\Delta')$ and $\Gamma=\pi(\Gamma')$ are fibers of~$\cX$. Thus, $\pi(s')$ is a union of some relations belonging to the set  $\{s\in S_{\Delta,\Gamma}:\ s\subseteq e_i\}$. From  Lemma~\ref{180421a}(2), it follows that  there is a unique~$s$ such that $\rho_{e'_i}(s)=\pi(s')$. Thus, $\pi(s')=s\in S$.\medskip

Let $k>1$. Take any $i\in\{1,\ldots,m\}$ for which $s'_i$ is irreflexive; in order to simplify notation, let $i=1$. 
Put
$$
x=s'_1\otimes 1_{\Omega_-(s'')}\qaq y=1_{\Omega_+(s'_1)}\otimes s'',
$$
where  $s''=s'_2\otimes\cdots\otimes s'_m$. Then $s'=x\cdot y$ by formula~\eqref{150421a}. Furthermore, $\pi(x),\pi(y)\in S$ by induction. Thus,
$$
\pi(s')=\pi(x\cdot y)=\pi(x)\cdot\pi(y).
$$
On the other hand, $x\subseteq f_1$ and $y\subseteq f'_1$. Consequently, $\pi(x)\subseteq e_1$ and $\pi(y)\subseteq e'_1$. Since $e^{}_1\perp e'_1$, we have $\pi(s')=\pi(x)\cdot \pi(y)\in S$, as required.\eprf

\crllrl{140121d}
Let $P=\{e_1,\ldots,e_m\}$ be an atomic  Cartesian decomposition of  a coherent configuration $\cX$. Assume that   $e$ and $f$ are parabolics of~$\cX$, and also $e\subseteq e_i$ and $f\subseteq e_j$ for some $1\le i\ne j\le m$. Then $e$ and $f$ strongly commute.
\ecrllr
\proof By Theorem~\ref{130121c}, we may assume that $\cX=\cX_1\otimes\cdots\otimes\cX_m$ and $P$ is the standard atomic  Cartesian decomposition of  $\Omega=\Omega_1\otimes\cdots\otimes\Omega_m$, where $\cX_i$ is a coherent configuration on~$\Omega_i$, $i=1,\ldots,m$. Then the rest of the proof is as in  Theorem~\ref{250121b}.\eprf

\thrml{280121c}
Let $\cX=\cX_1\otimes\cdots\otimes\cX_m$ and $P$ the standard atomic  Cartesian decomposition of $\cX$.  If the coherent configuration $\cX_i$ is indecomposable, $i=1,\ldots,m$, then $P$ is a $\preceq$-maximal atomic Cartesian decomposition of~$\cX$.
\ethrm
\proof  Assume on the contrary that $P$ is not $\preceq$-maximal. Then there exists an atomic Cartesian decomposition $R\ne P$ of~$\cX$, such that $P\preceq R$. By Lemma~\ref{040421c}(2), there exists a partition $\cR$ of the set~$M$ for which $P=R^\cR$.  Since $R\ne P$, one can find a non-singleton set $I\subseteq \{1,\ldots,|R|\}$ such that $R_I=e_i\in P$ for some $i\in M$. Then $\pr_I(R)$ is a nontrivial atomic Cartesian decomposition of $\cX_i$. By Theorem~\ref{130121c}, this implies that $\cX_i$ is decomposable, a contradiction.\eprf

\section{The uniqueness of maximal Cartesian decomposition}\label{270121a}

\subsection{Irredundant relations.}
Let $\cX=(\Omega,S)$ be a coherent configuration.  A  relation $s\in S$ is said to be {\it redundant} if  there exist irreflexive relations $x,y\in S$ such that 
\qtnl{120421a}
s=x\cdot y\qaq \grp{x}\wedge\grp{y}= 1_\Omega;
\eqtn
otherwise $s$ is said to be {\it irredundant}. Thus every redundant relation is the dot product of at least two irreflexive basis relations. One can also see that $s^*$ is redundant if and only if so is $s$.

\lmml{140121a2}
If $s\in S$ satisfies conditions~\eqref{120421a} for some $x,y\in S$, then $n_s=n_xn_y$.
\elmm
\proof We claim that $c_{xy}^s=1$. Indeed,  otherwise $ c_{xy}^s\ge 2$ and 
$$
|\alpha x\cap\beta y^*|=c_{xy}^s\ge 2,
$$ 
where $\alpha$ and $\beta$ are points such that $r(\alpha,\beta)=s$. For any two distinct points $\gamma$ and~$\delta$ of $\alpha x\cap\beta y^*$, we have
$$ 
r(\gamma,\delta)\subseteq (r(\gamma,\alpha)\cdot r(\alpha,\delta))\cap  (r(\gamma,\beta)\cdot r(\beta,\delta))=
(x^*\cdot x)\cap (y\cdot y^*)\subseteq\grp{x}\cap\grp{y}= 1_\Omega.
$$  
Consequently, $\gamma=\delta$, a contradiction. The claim is proved.\medskip

Next, since $s\in S$ and $s=x\cdot y$, we have $c_{xy}^t=0$ for all $t\ne s$. By the well-known identity for the intersection numbers (see, e.g., \cite[formula~(2.1.8)]{CP}),  this implies that
$$
n_s=n_sc_{xy}^s=\sum_{t\in S}n_tc_{xy}^t=n_xn_y,
$$ 
as required.\eprf\medskip

Given $s\in S$, put $d(s)=n_{s^{}}n_{s^*}$. Clearly, $d(s)=1$ if and only if the relation $s$ is thin. The invariant $d(s)$ is especially useful if the relation $s$ is redundant.

\crllrl{150421g}
Under the condition of Lemma~\ref{140121a2}, we have $d(s)=d(x)d(y)$.
\ecrllr
\proof  By Lemma~\ref{140121a2}, $n_s=n_xn_y$ and similarly $n_{s^*}=n_{x^*}n_{y^*}$. Thus,
$$
d(s)=n_{s^{}}n_{s^*}=(n_xn_y)\,(n_{x^*}n_{y^*})=(n_xn_{x^*})\, (n_yn_{y^*})=d(x)d(y),
$$
as required.\eprf\medskip

Not every coherent configuration has at least one irredundant relation: take a coherent configuration associated with regular elementary abelian group of rank at least~$2$. However, as the following lemma shows, every thick coherent configuration  has irredundant relations.

\lmml{140121a1}
Let $\cX$ be a thick coherent configuration, and let $d=\min_s d(s)$, where~$s$ runs over the irreflexive basis relations of~$\cX$. Then every $s\in S$ with $d(s)\le d$  is irredundant.
\elmm
\proof Let $s\in S$ be  such that $d(s)\le d$. Assume on the contrary that $s$ is redundant. Then there are irreflexive relations $x,y\in S$ such that formula~\eqref{120421a} holds. By Corollary~\ref{150421g}, we have 
$$
d\ge d(s)=d(x)d(y)\ge d^2,
$$ 
whence $d=1$. Consequently, $d(x)=d(y)=1$. Thus the relations $x$ and $y$ are thin. Since the coherent configuration $\cX$ is thick, they are reflexive, a contradiction.\eprf

\thrml{150421q}
Let $\cX$ be a thick coherent configuration and $s\in S$. Then there exist irredundant relations $s_1,\ldots,s_k\in S$ such that $s=s_1\cdot\cdots\cdot s_k$.
\ethrm
\proof By Lemma~\ref{140121a1}, we may assume that $d(s)>d$. If $s$ is irredundant, then we are done with $k=1$.  Otherwise  formula~\eqref{120421a} holds for some irreflexive $x,y\in S$. Then $d(s)=d(x)d(y)$ by Corollary~\ref{150421g}.  Since $\cX$ is thick and $x,y$ are irreflexive,  we have $d(x)>1$ and $d(y)>1$. Therefore, $d(x)<d(s)$ and $d(y)<s$. Now, the statement follows by induction.\eprf

\subsection{Atomic Cartesian decomposition $P^*$.} Let $\cX$ be a coherent configuration. We define a special set $P^*\subseteq E$ as the output of the following procedure.\medskip
 
\centerline{\bf Algorithm A}\medskip

\noindent{\bf Step 1.} Set $Q=\{\grp{s}:\ s\in S$ is irredundant$\}$.\smallskip
 
\noindent{\bf Step 2.} While there are distinct $e,f\in Q$ that are not strongly commute or  $e\not\perp f$, set 
$$
Q:=(Q\setminus\{e,f\})\cup \{e\vee f\}.
$$
\noindent{\bf Step 3.}  Output $P^*:=Q\setminus \{1_\Omega\}$.\medskip
 
Clearly, $P^*=\varnothing$ if and only if~$\cX$ contains no irredundant basis relations.  Note that, in general, $P^*$ may depend on the choice of~$e$ and~$f$ at each iteration of Step~$2$. In any case, the cardinality of the set~$Q$ defined  at Step 1 is at most $|S|$, and is reduced by one at each iteration of Step~2. Thus after at most $|S|$ iterations, Algorithm A arrives at Step~3.\medskip

The  lemma below immediately follows from Lemma~\ref{140121a1} and the description of Algorithm A.
 
\lmml{120421c}
Let $\cX$ be a thick coherent configuration. Then $P^*$ is a nonempty subset of~$E$,  satisfying the following conditions:
\nmrt
\tm{P0} $1_\Omega\not\in P^*$,
\tm{P1}  each irredundant relation of $\cX$ is contained in some parabolic from~$P^*$, 
\tm{P2} the parabolics from $P^*$ are pairwise strongly commute,
\tm{P3} $e\perp f$ for all distinct $e,f\in P^*$.
\enmrt
\elmm

In some special cases (e.g., if $\cX$ is the scheme of conjugacy classes of an arbitrary group, see~\cite[Example 2.4.3]{CP}), $P^*$  is an atomic  Cartesian decomposition of~$\cX$. It is not clear whether this is always true  (if it would be so, then the proof of the main theorems can essentially be simplified): indeed,  condition~(P3) is much weaker than that in Definition~\ref{229421f}. However, the following theorem implies, in particular, that~$P^*$ is an atomic Cartesian decomposition of~$\Omega$.

\thrml{170121a}
Let $\cX$ be  a thick coherent configuration on $\Omega$, and let $P^*\subseteq E$ be the output of the Algorithm~A.  Then
\nmrt
\tm{i} $P^*$  is an atomic Cartesian decomposition  of~$\Omega$,
\tm{ii} if $P$ is  an atomic  Cartesian decomposition  of~$\cX$, then $P\preceq P^*$.
\enmrt
\ethrm
\proof (i) By statements (P0) and (P2) of Lemma~\ref{120421c}, we need to verify~formula~\eqref{180421v} only. The first part of it follows from the lemma below for $Q=P^*$.

\lmml{220121a}
Let $Q$ be an arbitrary set of pairwise commuting parabolics of the coherent configuration~$\cX$. Assume that $e\perp f$ for all distinct $e,f\in Q$. Then $e\wedge e'=1_\Omega$ for all $e\in Q$, where~$e'$ is the $Q$-complement of $e$.
\elmm
\proof  Without loss of generality we may assume that $|Q|\ge 2$. The proof involves the induction on the cardinality of $Q$. For $|Q|=2$, the statement follows from Lemma~\ref{230421i}. Now let $|Q|\ge 3$ and $e\in Q$. \medskip

Assume on the contrary that there is an irreflexive basis relation $s\subseteq e\wedge  e'$. Then $d(s)>1$ by the thickness of~$\cX$. Without loss of generality, we may assume that $n_s\ge 2$; otherwise, we replace $s$ with~$s^*$.\medskip

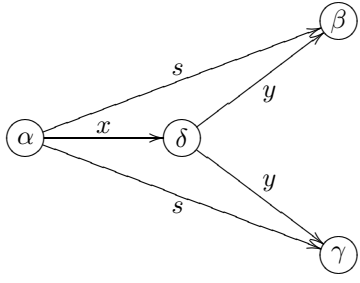
\begin{figure}[t]
$\xymatrix@R=30pt@C=15pt@M=0pt@L=2pt{
& & & & & & \VRT{\beta}\\
\VRT{\alpha} \ar[rrr]^*{x} \ar[rrrrrru]^*{s} \ar[rrrrrrd]_*{s} & & & \VRT{\delta} \ar[rrru]_*{y}  \ar[rrrd]^*{y} & & & \\
& & & & & & \VRT{\gamma}\\
}$
\caption{Configuration of the points $\alpha,\beta,\gamma,\delta$.}\label{fig1}
\end{figure}

Let $f\in Q\setminus\{e\}$ and $g$ the $Q$-complement of~$e\vee f$; in particular, $e'=f\vee g$. By the  induction applied to  $Q\setminus\{f\}$, we have 
\qtnl{220521a}
e\wedge g=1_\Omega.
\eqtn
Moreover, since $s\subseteq e'=f\vee g=f\cdot g$, there are basis relations $x\subseteq f$ and $y\subseteq g$ such that $s\subseteq x\cdot y$. It follows that $y\subseteq x^*\cdot s$.  Since also  $x^*\subseteq f$, $s\subseteq e$, and $f\perp e$, we have 
\qtnl{100521a}
y=x^*\cdot s.
\eqtn 
Take arbitrary $\alpha\in\Omega_-(s)$ and any two different points $\beta,\gamma\in \alpha s$ (recall that $n_s\ge 2$). Equality~\eqref{100521a} yields $\Omega_+(x^*)=\Omega_-(s)$, and hence there is a point $\delta\in \alpha x$ such that  the chosen points  form the configurations in Fig.~\ref{fig1}.
In particular, by~\eqref{220521a}, we obtain
$$
r(\beta,\gamma)\subseteq (s^*\cdot s) \cap (y^*\cdot y)\subseteq e\wedge g=1_\Omega,
$$
It follows that $\beta=\gamma$, a contradiction.\eprf\medskip

To complete the proof of (i), assume on the contrary that $e\vee e'\ne\bone_\Omega$.  Then there is $s\in S$ such that $s\cap (e\vee e')=\varnothing$. Note that $s$ is redundant, because by the condition~(P1) every irredundant relation is contained in either~$e$ or $e'$. By Theorem~\ref{150421q}, $s=x_1\cdot  x_2\cdots x_k$ for some $k\ge 2$ and irredundant $x_1,\ldots,x_k\in S$.  Again by condition~(P1),  $x_i\subseteq e\vee e'$ for all~$i$. Thus, 
$$
s=x_1\cdot  x_2 \cdots  x_k\subseteq  (e\vee e')\cdot   (e\vee e')\cdots   (e\vee e')= e\vee e',
$$
a contradiction.\medskip

(ii)\quad  Let $P=\{e_1,\ldots,e_m\}$ be an atomic Cartesian decomposition  of~$\cX$. It suffices to verify that  if $Q\subseteq E$ is  obtained at a certain step of Algorithm~A, then $P\preceq Q$, i.e., every $e\in Q$ is contained in some~$e_i\in P$.\medskip

Let $Q$ be constructed at Step~1. Assume on the contrary that $e$ is  contained in no $e_i\in P$. Denote by $J$ a maximal subset of $M$ such that $e\not\subseteq P_J$. Then $1\le |J|<m$. It follows that $e\subseteq e_i\cdot P_J$
for every $i\in M\setminus J$.  On the other hand, $e=\grp{s}$ for some irredundant $s\in S$. By formula~\eqref{070221a}, we have $e_i\perp P_J$, because $P_J\subseteq e_i'$ and $e^{}_i\perp e'_i$. Consequently,
$$
s=x\cdot y,\quad x\subseteq e^{}_i,\ y\subseteq P_J.
$$
Moreover, none of $x$, $y$ is reflexive: if $x$ is reflexive, then $s=y\subseteq P_J$ and hence $e=\grp{s}\subseteq P_J$, in contrast to the choice of~$J$, whereas if $y$ is reflexive, then $s=x\subseteq e_i$ and hence $e=\grp{s}\subseteq e_i$, in contrast to the assumption. Finally,
$$
\grp{x}\wedge\grp{y}\subseteq e^{}_i\wedge P_J\subseteq  e^{}_i\wedge e'_i=1_\Omega.
$$
Consequently, the relation $s$ is redundant, a contradiction.\medskip

Now let $Q$ be obtained at a certain iteration of Step~2, not just after Step~1. Denote by $R$ the set obtained at the previous iteration. Then $Q=(R\setminus\{f, g\})\cup\{f\vee g\}$  for some distinct  $f, g\in R$ which are not strongly commute  or  $f\not\perp g$. By induction, we may assume that $P\preceq R$. Consequently, $e_i\supseteq f$ and $e_j\supseteq g$ for some $i,j\in M$.\medskip

Let $e$ be an arbitrary equivalence relation of~$Q$. If  $e\ne f\vee g$, then $e\in R$ and hence $e$ is contained in some element of~$P$, as required. Let  $e=f\vee g$. If, in addition, $i=j$, then $e=e_i\in P$. Now let $i\ne j$. Then  $e_j\subseteq e'_i$. Since $e^{}_i\perp e'_i$, formula~\eqref{070221a} shows that $f\perp g$. By the choice of $f$ and $g$, this means that  $f$ and $g$ are not strongly commute. However, this contradicts Corollary~\ref{140121d}.\eprf\medskip

The following sufficient condition  for a thick coherent configuration to be indecomposable is an immediate consequence of Theorem~\ref{170121a}(2).

\crllrl{250121a}
A thick coherent configuration $\cX$ is indecomposable if $|P^*|=1$.
\ecrllr

\subsection{Proofs of Theorem~\ref{080619a} and Corollary~\ref{110421c}.}\phantom{x}\medskip

{\bf Proof of Theorem~\ref{080619a}}. Because of the poset anti-isomorphism between the Cartesian and atomic Cartesian decompositions of~$\cX$, we prove the statement for the latter ones.  Clearly, $\cX$ has at least one  atomic Cartesian decomposition, the trivial one. Assume that $P$ and $Q$ are  $\preceq$-maximal atomic Cartesian decompositions of~$\cX$.  By  Theorem~\ref{170121a}, the set $P^*$ is  an atomic Cartesian decomposition of~$\Omega$, and
\qtnl{110421q}
P\cap Q\preceq P^*.
\eqtn
Consequently, $P\cap Q$ is also an atomic Cartesian decomposition of~$\Omega$ (Lemma~\ref{260121a}). It remains to verify the claim below, because then $P=P\cap Q=Q$ by the maximality of $P\preceq P\cap Q$ and $Q\preceq P\cap Q$. \medskip

{\bf Claim.} {\it  $P\cap Q$ is an atomic Cartesian decomposition of~$\cX$.}\medskip

\proof Let $e\in P\cap Q$ and $e'$ the $(P\cap Q)$-complement of~$e$. There are unique $f\in P$ and $g\in Q$ such that $e=f\wedge g$. Denote by $f'$ the $P$-complement of~$f$. Since $P\cap Q$ is  an atomic Cartesian decomposition of~$\Omega$, Lemma~\ref{040421c}(2) implies that $f,f',g\in L(P\cap Q)$ and also
$f=e\cdot\ov f$ for some $\ov f\in L(P\cap Q)$. Thus,
\qtnl{140421a}
e'=\ov f\cdot f'.
\eqtn

Let  $x,y\in S$ be arbitrary relations such that $x\subseteq e$ and $y\subseteq e'$. Note that $\ov f\perp f'$, because $\ov f\subseteq f$ and $f\perp f'$. Since the parabolics $\ov f,f'\in L(P\cap Q)$ commute and in view of~\eqref{140421a}, this implies that there exist $\ov y,y'\in S$ such that $\ov y\subseteq\ov f$, $y'\subseteq f'$, and 
$$
y=\ov y\cdot y'.
$$

In order to verify $e\perp e'$, we prove that $x\cdot y\in S$ provided that $x\cdot y\ne\varnothing$.  We have $x\subseteq e\subseteq g$ and $\ov y\subseteq \ov f\subseteq g'$. Since $g\perp g'$, this implies that $x\cdot \ov y\in S$ (note that $x\cdot \ov y\ne\varnothing$, for otherwise $x\cdot y=x\cdot \ov y\cdot y'=\varnothing$). Taking into account that $x\cdot \ov y\subseteq f$, $y'\subseteq f'$, and $f\perp f'$, we obtain
$$
x\cdot y=x\cdot (\ov y\cdot y')=(x\cdot \ov y)\cdot y'\in S,
$$
as required.\medskip

Let us verify that $e$ and $e'$ strongly commute. Note that $x$ is contained in both $g$ and $f$. Since $g$ strongly commutes with $g'\supseteq \ov f$ and $f$ strongly commutes with $f'\supseteq y'$, there are basis relations $\ov z\subseteq g'$ and $z'\subseteq f'$ such that
$$
x\cdot y=x\cdot (\ov y\cdot y')=(x\cdot \ov y)\cdot y'=\ov z\cdot x\cdot y'=\ov z\cdot z'\cdot x.
$$ 
Take an arbitrary $z\in S$ contained in $\ov z\cdot z'\subseteq g'\cdot f'\subseteq e'$ and such that $z\cdot x=\varnothing$ if $x\cdot y=\varnothing$, and $z\cdot x$ intersects $x\cdot y$ otherwise. In the former case, $z\cdot x=x\cdot y$. In the latter case, we note that $e'\perp e$ by the above, and hence $x\cdot y, z\cdot x\in S$. Thus, $z\cdot x=x\cdot y$ and hence $x\cdot e'=e'\cdot x$. Similarly, one can verify that $e\cdot y=y\cdot e$ \eprf\medskip

{\bf Proof of Corollary~\ref{110421c}.} Let $\cX$ be a thick coherent configuration. Without loss of generality, we may assume that $\cX$ is decomposable. Suppose that  we are given coherent configurations 
$$
\cX_1\otimes\cdots \otimes \cX_l\cong\cX\qaq \cY_1\otimes\cdots \otimes \cY_m\cong\cX
$$ 
the factors of which are  indecomposable of degree greater than~$1$. The assumption implies that $l>1$ and $m>1$. By Theorem~\ref{250121b}, there are atomic Cartesian decompositions $P=\{e_1,\ldots,e_l\}$ and $Q=\{f_1,\ldots,f_m\}$ of the coherent configuration~$\cX$, such that
\qtnl{270121q}
\cX_i\cong\cX_{\Omega/e'_i}\qaq \cY_j\cong\cX_{\Omega/f'_j} 
\eqtn
for $i=1,\ldots,l$ and $j=1,\ldots,m$. Moreover, these decompositions are $\preceq$-maximal by Theorem~\ref{280121c}. By Theorem~\ref{080619a}, this yields $P=Q$. In particular, $l=m$, and, after a suitable enumeration, $e_i=f_i$ for $i=1,\ldots,l$. In view of~\eqref{270121q}, this implies that  $\cX_i\cong\cY_i$ for all~$i$.\eprf

 \section{Finding the maximal Cartesian decomposition}\label{110421v2}
 
In this section we prove Theorem~\ref{080619a1} by constructing a polynomial-time algorithm to find the maximal Cartesian decomposition of a thick coherent configuration~$\cX$. At the input of the algorithm, $\cX$ is given by a list of basis relations and each of which is considered as a directed graph. In this representation, for any two relations $r$ and $s$ of $\cX$, one can efficiently find  $\grp{r}$, $r\cap s$, $r\cup s$, and $r\cdot s$, and check whether $r$ and $s$ commute; the corresponding algorithms are standard and can be implemented with the help of basic procedures described in, e.g.,~\cite[Chap.~VI]{CLRS2009}. 

\lmml{060221a}
Let $\cX$ be a coherent configuration of degree $n$. Then  in polynomial time in~$n$, one can test
\nmrt
\tm{i} whether a given relation $s\in S$ is irredundant,
\tm{ii} whether given parabolics $e,f\in E$ strongly commute, or whether $e\perp f$.
\enmrt
\elmm
\proof To test (i), one need to verify that $s=x\cdot y$ for no irreflexive $x,y\in S$ such that $\grp{x}\cap\grp{y}=1_\Omega$. By the above remarks, this  can be done efficiently for a fixed pair~$x$ and~$y$. Since there is at most~$|S|^2\le n^4$ of such pairs, we are done. Testing~(ii) requires verifying equalities of the form $r\cdot s=s'\cdot r'$ for fixed $r, s$ and all $r',s'\in S$, or  relations \eqref{060221c};  both  can be done efficiently by the above remarks. \eprf\medskip
 
\crllrl{060221b}
In the condition of Lemma~\ref{060221a}, the set $P^*$ defined by Algorithm~A can be found in polynomial time in~$n$.
\ecrllr

A {\it decomposability certificate} of $\cX$ is an atomic Cartesian decomposition $P$ of $\cX$ such that $|P|\le 2$ and $|P|=2$ if and only if $\cX$ is decomposable. The construction of such a certificate is sufficient to find the $\preceq$-maximal atomic Cartesian decomposition of $\cX$. The input coherent configuration in the following algorithm is assumed to be thick; for nonthick coherent configurations the algorithm is not always correct.\medskip

\centerline{\bf Algorithm B}\medskip
 
{\bf Input:} a thick coherent configuration $\cX$ on $\Omega$.
 
{\bf Output:} a  decomposability certificate $P$ of $\cX$.\medskip
 
\noindent{\bf Step 1.}  Find  the set $P^*\subseteq E$ with the help of Algorithm A; put $M^*=M(P^*)$.\smallskip

\noindent{\bf Step 2.}  For each  proper subset $I\subset M^*$, test whether $P:=\{P^*_I,P^*_{M^*\setminus I}\}$ is an atomic Cartesian decomposition of $\cX$, and if so, then output~$P$.\smallskip
 
\noindent{\bf Step 3.} Output $P:=\{\bone_\Omega\}$.\medskip

\lmml{090421c}
Algorithm B constructs a decomposability certificate $P$ correctly and in time polynomial in~$n$. 
\elmm
\proof Let us prove the correctness.  If $P$ is constructed at Step~2, then $|P|=2$ and $\cX$ is decomposable by Theorem~\ref{130121c}. Let $P$ be constructed at Step~3. Assume on the contrary that there is a nontrivial atomic Cartesian decomposition~$R$ of~$\cX$. By Theorem~\ref{170121a}, the set $P^*$ found  at Step~1 is an atomic Cartesian decomposition of $\Omega$ and  $R\preceq P^*$.  By Lemma~\ref{040421c}(2)  applied to $P=P^*$,  there is a partition $\cR$ of~$M^*$ such that 
$$
R=(P^*)^\cR.
$$ 
Now, let $I\in\cR$. Then the above equality implies that $P^*_I\in R$. Consequently, the set $P$ defined at Step~$2$ is an atomic  Cartesian decomposition of~$\cX$. This means that Algorithm~B stops at Step~2, a contradiction.\medskip

Let us estimate the running time of Algorithm B on an input of degree~$n$. We note that Step~1 can be implemented in polynomial time by Corollary~\ref{060221b}. At Step~2, one needs to inspect all proper subsets of~$M^*$. Since $P^*$ is an atomic  Cartesian decomposition of~$\Omega$, formula~\eqref{080221c} implies that 
$$
n=|\Omega|=\prod_{i\in M^*}|\Omega/e_i'|\ge 2^{|M^*|}.
$$
Assuming $n>1$, we have  $|M^*|\le \log n$ and hence the number of proper subsets of~$M^*$ is at most~$n$. Moreover, testing that for a given $I\subset M^*$ the set $\{P^*_I,P^*_{M^*\setminus I}\}$ is an atomic Cartesian decomposition of~$\cX$ can be done in polynomial time by Lemma~\ref{060221a}(ii). This proves that Step~2 and then the entire algorithm is polynomial-time.\eprf\medskip

{\bf Proof of Theorem~\ref{080619a1}.} Recall that a tensor decomposition of  coherent configuration $\cX$ on~$\Omega$ into indecomposable components is a set $C$ of indecomposable coherent configurations $\cX_i$ on $\Omega_i$, $i=1,\ldots,m$, and an isomorphism $\pi:\Omega\to\Omega_1\times\cdots\times\Omega_m$ from $\cX$ to~$\cX_1\otimes\cdots\otimes\cX_m$.\medskip

\centerline{\bf Algorithm C}\medskip

{\bf Input:} a thick coherent configuration $\cX$ on $\Omega$.

{\bf Output:} a   tensor decomposition $(C,\pi)$ of $\cX$ into indecomposable components.\medskip

\noindent{\bf Step 1.}  Find  a decomposability certificate $P$ of $\cX$ with the help of Algorithm B.\smallskip

\noindent{\bf Step 2.}  If $|P|=1$, then output $C=\{\cX\}$ and $\pi=\id_\Omega$. Otherwise, let $P=\{e_1,e_2\}$, $\cX_1=\cX_{\Omega/e_1}$ and $\cX_2=\cX_{\Omega/e_2}$.\smallskip

\noindent{\bf Step 3.} For $i=1,2$, recursively find a   tensor decomposition $(C_i,\pi_i)$ of $\cX_i$  into indecomposable components: $C_i=\{\cX_{i,1_{}},\ldots,\cX_{i,m_i}\}$ and $\pi_i$ is an isomorphism from~$\cX_i$ to  $\cX_{i,1}\otimes\cdots\otimes\cX_{i,m_i}$.\smallskip

\noindent{\bf Step 4.} Set $\pi'$  to be the isomorphism from $\cX_1\otimes\cX_2$ to $\cX_{1,1}\otimes\cdots\otimes\cX_{1,m_1}\otimes\cX_{2,1}\otimes\cdots\otimes\cX_{2,m_2}$, induced by the isomorphisms $\pi_1$ and $\pi_2$.\smallskip

\noindent{\bf Step 5.} Output $C=C_1\cup C_2$ and  the composition $\pi$  of $\pi_P$ with $\pi'$.\medskip

The correctness of the output at Step~2 follows from Lemma~\ref{090421c}. Furthermore,  the recursive calls at Step~3 are correct. Indeed, if one of the coherent configurations~$\cX_1$, $\cX_2$ is not thick, say $\cX_1$ contains 
an irreflexive thin basis relation~$s_1$, then~$\cX$ is not thick: $S$ has the irreflexive thin relation $\pi_P^{-1}(s_1\otimes s_2)$, where~$s_2$ is an arbitrary reflexive relation of~$\cX_2$. Thus the coherent configurations $X_{i,j}$ found at Step~3 are indecomposable by induction. It remains to note that the bijection $\pi$ constructed at Step~5 is a composition of two isomorphisms (the fact that $\pi_P$ is an isomorphism follows from Theorem~\ref{130121c}).\medskip

To prove that Algorithm C is polynomial-time, it suffices to bound from above the number $c(n)$ of the recursive calls at Step~3 by a polynomial in~$n$: indeed, each  step except for Step~1 can obviously be implemented efficiently, whereas for Step~1, this is true by Lemma~\ref{090421c}. Finally, at Step~3, we have $|\Omega|=|\Omega_1\times\Omega_2|$ and hence $n=n_1\,n_2$, where $n_1=|\Omega_1|$ and $n_2=|\Omega_2|$. Thus,
$$
c(n)=2+c(n_1)+c(n_2).
$$
This immediately implies that $c(n)=\log n$, as required.\eprf

\end{document}